\documentclass[12pt]{amsart}

\def\CB{\mathcal {B}}
\def\CL{\mathcal {L}}
\def\CD{\mathcal {D}}
\def\H{\mathbb{H}}
\def\C{\mathbb{C}}
\def\N{\mathbb{N}}

\def\R{\mathbb{R}}

 \newtheorem{thm}{Theorem}[section]

\newcommand{\be}{\begin{equation}}
\newcommand{\ee}{\end{equation}}
\newcommand{\bea}{\begin{eqnarray}}

\newcommand{\eea}{\end{eqnarray}}
\newcommand{\Bea}{\begin{eqnarray*}}
\newcommand{\Eea}{\end{eqnarray*}}

\newcounter{cnt1}
\newcounter{cnt2}
\newcounter{cnt3}
\newcommand{\blr}{\begin{list}{$($\roman{cnt1}$)$}
 {\usecounter{cnt1} \setlength{\topsep}{0pt}
 \setlength{\itemsep}{0pt}}}
\newcommand{\bla}{\begin{list}{$($\alph{cnt2}$)$}
 {\usecounter{cnt2} \setlength{\topsep}{0pt}
 \setlength{\itemsep}{0pt}}}
\newcommand{\bln}{\begin{list}{$($\arabic{cnt3}$)$}
 {\usecounter{cnt3} \setlength{\topsep}{0pt}
 \setlength{\itemsep}{0pt}}}
\newcommand{\el}{\end{list}}

\date{}
\begin{document}

\title[Paley-Wiener theorem ]
{A Paley-Wiener theorem for the inverse  \\
\vskip .5em Fourier transform on some\\
\vskip .5em  homogeneous spaces \\
\vskip 1.5em {\tt By} }

\author[Thangavelu]{S.\ Thangavelu}
\address{Department of Mathematics,\\
Indian Institute of Science,\\
Bangalore 560 012, India {\it E-mail~:} {\tt
veluma@math.iisc.ernet.in}}

\maketitle

\begin{abstract}
We formulate and prove  a version of Paley-Wiener theorem for the inverse 
Fourier transforms on noncompact Riemannian symmetric spaces and Heisenberg 
groups. The main ingredient in the proof is the Gutzmer's formula.
\end{abstract}

\section{Introduction}
\setcounter{equation}{0}

The classical Paley-Wiener theorem for the Euclidean Fourier transform
characterises compactly supported functions on $ \R^n $ in terms of holomorphic
properties of their Fourier transforms. Analogues of Paley-Wiener theorem have
been proved in the context of Fourier transforms on Lie groups. One such 
result is a theorem of Gangolli [5] for the spherical Fourier transform on
noncompact Riemannian symmetric spaces. However, there are no satisfactory
results available in certain cases. One such example is the case of the
Fourier transform on the Heisenberg group $ \H^n.$ Here the Fourier transform
is operator valued, parametrised by non-zero real numbers. When a function $ f
$ on $ \H^n$ is compactly supported it is not possible to extend the Fourier
transform $ \hat{f}(\lambda) $ as an operator valued entire function. There
are some versions of Paley-Wiener theorem for the Heisenberg group which
treat the central and non-central variables separately, see [1] and [12]. The
situation of general nilpotent Lie groups is much more difficult.

In 2000, Pasquale [10] considered the problem of characterising
functions on a non-compact symmetric space $ X = G/K $ whose spherical
Fourier transforms are compactly supported. When $ G $ is a complex semisimple
Lie group or of rank one she showed that $ K-$ biinvariant functions whose 
spherical Fourier transforms
are compactly supported can be extended to the complexification of $ X $ as
meromorphic functions leading to a Paley-Wiener theorem for the inverse
spherical Fourier transform. One of the main results of this paper is such
a theorem for Helgason Fourier transform of general functions on $ X .$ 
The main ingredient in the proof
(which also motivates the formulation) is Gutzmer's formula proved by Faraut
[3] for functions holomorphic in a domain, called the complex crown, contained
in the complexification of $ X.$ Our result is similar in spirit to the
characterisation of the image of the heat kernel transform studied by
Kroetz et al [9].

Instead of treating $ \H^n $ as a nilpotent Lie group we may consider it
as  a homogeneous space of a bigger group, namely the Heisenberg motion group
$ G_n $ which is the semidirect product of $ \H^n $ with the unitary group
$ U(n).$ Thus $ \H^n = G_n/U(n) $ and we treat functions on $ \H^n $ as right
$U(n)-$ invariant functions on $ G_n.$ With this view point the Fourier
transform of $ f $ on $ \H^n $ is considered as a functions of two
variables $ \lambda $ and $ k.$ Here $ \lambda $ is a nonzero real number
and $ k $ is a non-negative integer. For each such pair there is a unitary
representation of $ G_n $ denoted by $ \rho_k^\lambda $ and we consider
the operators $ \rho_k^\lambda(f) $ as parametrised by the point $ (\lambda,
(2k+n)|\lambda|) $ from the Heisenberg fan which is the spectrum of the
sublaplacian. Our Paley-Wiener theorem for the Heisenberg group characterises
functions for which $ \rho_k^\lambda(f) $ is supported in $ |\lambda| \leq a $
and $ (2k+n)|\lambda| \leq b^2.$ The proof requires an analogue of Gutzmer's
formula for the Heisenberg motion group which has been proved recently [15].

The theorem of Kroetz et al [9] and our Paley-Wiener theorem both involve
a certain peseudo-differential shift operator $ D $. As was shown in [9]
the operator is inevitable in characterising the image of the heat kernel
transform. When the group $ G $ is complex the operator $ D $ is simple
( multiplication by a Jacobian factor) but otherwise it is quite complicated.
Interestingly enough our Paley-Wiener theorem for the Heisenberg group also
involves a similar operator $ \CD.$ The operator $ D $ has the effect of
replacing the elementary spherical function $ \varphi_\lambda $ by the
Weyl symmetrised exponential $ \psi_\lambda.$ The same is true of $ \CD.$
It, in effect, changes the Laguerre functions $ \varphi_k^\lambda $ into
Bessel functions. Without these operators, certain orbital integrals are not
entire functions of exponetial type.

It is worthwhile to see how our version of Paley-Wiener theorem looks like 
for the Euclidean Fourier transform. Let $ f $ be a Schwartz class function
on $ \R^n $ and consider the Fourier transform $ \hat{f}.$ When $ \hat{f} $
is supported in $ |\xi| \leq a $ $ f $ extends to $ \C^n $ as an entire 
function.  Let $ G = M(n) $ be the Euclidean motion group acting on $ \R^n$
which has a natural extension to $ \C^n.$ Then it is easy to see that
the following Gutzmer's formula is valid:
$$ \int_G |f(g.z)|^2 dg = c_n \int_0^\infty \int_{S^{n-1}}|\hat{f}(\lambda 
\omega)|^2 \varphi_\lambda(2iy) \lambda^{n-1}d\omega d\lambda $$
where $ \varphi_\lambda(iy) = (\lambda |y|)^{-\frac{n}{2}+1}
J_{\frac{n}{2}-1}(i\lambda |y|)$ and $ z = x+iy.$ From the above it is clear
that the orbital integral $ \int_G |f(g.z)|^2 dg $ satisfies the estimate
$$  \int_G |f(g.z)|^2 dg \leq C e^{2a|y|}.$$ Conversely, if a Schwartz
function $ f $ extends to $ \C^n$ as an entire function and the orbital
integral satisfies the above estimate then $ \hat{f} $ is supported in
$ |\xi| \leq a.$ This follows easily from the Gutzmer's formula.

The plan of the paper is as follows. In the next section we treat the
Helgason Fourier transform on non-compact Riemannian symmetric spaces.
In Section 3 we consider the inverse Fourier transform on the Heisenberg
group.

\section{Noncompact Riemannian symmetric spaces}
\setcounter{equation}{0}

In this section we formulate and prove a Paley-Wiener theorem for the inverse
Helgason Fourier transform on a Riemannian symmetric space of noncompact type.
We follow the standard notations; in fact we closely follow Kroetz et al
in setting up the notation and we refer to the same for any undefined term.
Let $ X = G/K $ be a homogeneous space where $ G $ is a semisimple Lie group
and $ K $ a maximal compact subgroup. Considering the Iwasawa decomposition
$ G = KAN $ we let $ M $ be the centraliser of $ A $ in $ K.$ We define
$ B = K/M $ and consider the Helgason Fourier transform 
$$ \hat{f}(\lambda,b) = \int_X f(x) e^{(-\lambda+\rho,A(x,b))} dx $$
 where $ \lambda \in i\bf{a}^* $  and $ b \in B.$ Here $ \bf{a} $ and $ A(x,b)
$ have the usual meaning. The inversion formula valid for suitable functions
reads as follows:
$$ f(x) = \int_{i\bf{a}^*} \left( \int_B  e^{(\lambda+\rho,A(x,b))}db\right)
|c(\lambda)|^{-2} d\lambda.$$

For every $ \lambda \in \bf{a}_\C^*$ ,$ b \in B $ the function $ x \rightarrow
e^{(\lambda+\rho,A(x,b))} $ has a holomorphic extension to a domain $ \Xi$
in the complexification $ X_\C. $ This domain, called the complex crown of
$ X $ is defined as follows. Let $ \bf{g} $ be the Lie algebra of $ G $ with
the Cartan decomposition $ \bf{g} = \bf{k}+\bf{p} .$  Let $ \bf{a} $ be a
Cartan subspace with $ \Sigma $ the associated system of restricted roots. The
complex crown $ \Xi $  is a $ G-$ invariant domain in $ X_\C = G_\C/K_\C $
defined by $ \Xi = G \exp(i\Omega).x_0, x_0 = eK $ where $ \Omega =
\{ H \in \bf{a}: |\alpha(H)| <\frac{\pi}{2}, \alpha \in \Sigma \}.$ (see 
Akhiezer-Gindikin [1] and Kroetz-Stanton [7]). It follows that a function $ f
\in L^2(X) $ whose Fourier transform $ \hat{f}(\lambda,b) $  has compact 
support admits a holomorphic extension to $ \Xi $ which is given by
$$ f(z) = \int_{i\bf{a}^*} \left( \int_B  e^{(\lambda+\rho,A(z,b))}db\right)
|c(\lambda)|^{-2} d\lambda.$$

Let $ \varphi_\lambda , \lambda \in \bf{a}^*$ be the spherical functions 
on $ G $ given by the integral
$$ \varphi_\lambda(g) = \int_B e^{(\lambda+\rho,A(x,b))}db .$$ In [7] Kroetz-
Stanton proved that for $\lambda \in i\bf{a}^* $ the function $ H \rightarrow
\varphi_\lambda(\exp iH) $ admits a holomorphic continuation to the tube
$ \bf{a}+2i\Omega.$ In order to state Gutzmer's formula and formulate a Paley
-Wiener theorem we need to recall the definition of orbital integrals developed
by Gindikin et al [6]. For a function $ h $ on $ \Xi $ suitably decreasing at
the boundary and $ Y \in 2\Omega $ we define
$$ O_h(iY) = \int_G h(g\exp(\frac{i}{2}Y).x_0) dg.$$
Let $ G(\Xi) $ be the space of all holomorphic functions on the complex
crown. Then in [3] Faraut has established the following formula known
as Gutzmer's formula.

\begin{thm}
Let $ f \in G(\Xi) $ be such that for all $ H \in \Omega $,
$$ \int_G |f(g\exp(iH).x_0)|^2 dg \leq M $$
for some constant $ M.$ Then for all $ H \in \Omega $ we have
$$\int_G |f(g\exp(iH).x_0)|^2 dg $$ 
$$ = \int_{i\bf{a}^*} \left( \int_B  |\hat{f}(\lambda,b)|^2 db\right) 
\varphi_\lambda(\exp(2iH)) |c(\lambda)|^{-2} d\lambda.$$
\end{thm}

In [9] Kroetz et al has used this Gutzmer's formula to characterise the
image of the heat kernel transform. Let $ k_t(x) $ stand for the heat
kernel associated to the Laplace-Beltrami operator on $ X $ which is
a $ K$-biinvariant function given by the integral
$$ k_t(x) = \int_{i\bf{a}^*} e^{-t(|\lambda|^2+|\rho|^2)}\varphi_\lambda(x)
 |c(\lambda)|^{-2} d\lambda.$$
It is clear that $ k_t $ has a holomorphic extension to the complex crown.
If $ f \in L^2(X) $ the function $ H_tf(x) = f*k_t(x) $ which solves the
heat equation with initial condition $ f $ also extends to $ \Xi $ as a 
holomorphic extension. Let $ im H_t $ stand for the image of the above
transform, $ f(x) \rightarrow H_tf(z) $ called the heat kernel transform.
For the Euclidean Laplacian the corresponding image turned out to be a
weighetd Bergman space; the same is true for compact symmetric spaces.
However, in [9] Kroetz et al proved that $ im H_t $ is not a weighted Bergman
space. Instead they obtained the following characterisation.

In order to state their result we need to set up some more notation. Let $ W $
be the Weyl group and consider the Weyl symmetrised exponential function
$$ \psi_\lambda(Z) = \sum_{w \in W} e^{(\lambda,wZ)},~~ Z \in \bf{a}_\C,~~
 \lambda
\in i\bf{a}^*.$$ If a holomorphic function $ h $ on the tube domain $ \bf{a}
+2i\Omega $ has the representation
$$ h(Z) = \int_{i\bf{a}^*} g(\lambda) \varphi_\lambda(\exp(Z).x_0)
|c(\lambda)|^{-2} d\lambda $$
then we define
$$ Dh(Z) = \int_{i\bf{a}^*} g(\lambda)\psi_\lambda(Z)
|c(\lambda)|^{-2} d\lambda .$$ Under some conditions on $ g $ this is well
defined, see [9]. It is known that $ D $ is a pseudo-differential shift 
operator which has a simpler form when the group $ G $ is complex. It can
be expressed in terms of Abel transform and Fourier multipliers. Using this
operator $ D $ the following characterisation of the heat kernel transform
was obtained in [9].

\begin{thm}
A function $ F \in G(\Xi) $ belongs to $ im H_t $ if and only if
$$ \int_{\bf{a}} DO_{|F|^2}(iY) w_t(Y) dY < \infty $$
where $w_t$ is given in terms of the Euclidean heat kernel as
$$ w_t(Y) = |W|^{-1}(2\pi t)^{-\frac{n}{2}}e^{2t|\rho|^2}e^{-\frac{|Y|^2}{2t}}
.$$
\end{thm}

This theorem is an easy consequence of the Gutzmer's formula. Our Paley-Wiener
theorem is very similar in spirit to the above theorem.

\begin{thm}
Let $ f $ be a function in  $ L^2(X).$ Then the Helgason Fourier transform 
$ \hat{f}(\lambda,b) $ is
supported in $ |\lambda| \leq R $ if and only if $ f $ has a holomorphic extension $ F \in G(\Xi)$ which satisfies the estimate
$$  DO_{|F|^2}(iY) \leq C e^{2R|Y|} $$
for some constant $ C $ independent of $ Y.$
\end{thm}

{\bf Proof}:  First assume that  $ \hat{f}(\lambda,b) $ is compactly supported
in  $ |\lambda| \leq R .$ From the inversion formula for the Helgason Fourier
transform it is clear that $ f $ can be holomorphically extended to $ \Xi.$
If $ F $ is the extension then by Plancherel theorem it follows that $ F \in
im H_t $ for all $ t >0.$ Moreover, Gutzmer's formula can be applied and we
get
$$  DO_{|F|^2}(iY) 
= \int_{i\bf{a}^*} \left(\int_B |\hat{f}(\lambda,b)|^2 db \right)
\psi_\lambda(2iY) |c(\lambda)|^{-2} d\lambda.$$
This gives the estimate, 
$$ DO_{|F|^2}(iY)
\leq C \|f\|_2^2 ~e^{2R|Y|} $$
as  $ \hat{f}(\lambda,b) $ is supported in $ |\lambda| \leq R $ and 
$|\psi_\lambda(iY)| \leq C e^{|\lambda||Y|}.$

Conversely, assume that $ F $ satisfies the hypothesis of the theorem. Then
it is easy to see that $ F \in im H_t $ for every $ t> 0. $ More
precisely, for every $ t > 0 $ we have
$$ \int_{\bf{a}}DO_{|F|^2}(iY)w_t(Y) dY \leq C e^{2t|\rho|^2}
P(R,t) e^{2 t R^2}$$
where $ P $ is some polynomial. Consider the integral
$$ \int_{|\lambda|\geq R+\epsilon} \left(\int_B |\hat{f}(\lambda,b)|^2 db 
\right) |c(\lambda)|^{-2} d\lambda$$
$$ \leq e^{-2t(R+\epsilon)^2}\int_{i\bf{a}^*} \left(\int_B |\hat{f}(\lambda,b)|^2 db \right)
e^{2t|\lambda|^2} |c(\lambda)|^{-2} d\lambda.$$ By the above and Gutzmer's
formula, we get the estimate
$$ \int_{|\lambda|\geq R+\epsilon} \left(\int_B |\hat{f}(\lambda,b)|^2 db
\right) |c(\lambda)|^{-2} d\lambda \leq C e^{-2t(R+\epsilon)^2}P(R,t)
e^{2tR^2}.$$
By letting $ t $ tend to infinity we conclude that $ \hat{f}(\lambda,b) $
vanishes almost everywhere for $ |\lambda| \geq R+\epsilon.$ As $ \epsilon
$ is arbitrary $ \hat{f}(\lambda,b) $ is supported in $ |\lambda| \leq R $
proving the theorem.

We conclude this section with the following remarks. For each $ t> 0 $ the
image $ im H_t $ is a Hilbert space with the norm
$$ \|F\|_t^2 = \int_{\bf{a}} DO_{|F|^2}(iY) w_t(Y) dy.$$ As shown in [9]
$  \|F\|_t = \|f\|_2 $ if $ F = f*k_t.$ Let $ \Delta $ be the 
Laplace-Beltrami operator, taken to be non-negative so that $ e^{-t \Delta}
f = f*k_t.$ Let us define $ \bf{H} $ to be the intersection of all
$ im H_t , t>0.$ If $ L^2_b(X) $ stand for the subspace of $ L^2(X) $ with
compactly supported Helgason Fourier transforms then it is clear that 
$f $  is the restriction of an  $ F \in \bf{H}.$ The above theorem can be 
viewed as one charactersing the image
of $ L^2_b(X) $ under the heat kernel transform.

\section{Fourier transform on the Heisenberg group}
\setcounter{equation}{0}

In this section we consider the Heisenberg group as the homogeneous space
$ G_n/U(n) $ where $ G_n $ is the Heisenberg motion group. The general
references for this section are the papers Kroetz et al [8] and [15]. See also
the monographs [4] and [13]. We take $ \H^n $ to be $ \C^n \times \R $ with
group law $ (z,t)(w,s) = (z+w,t+s+\frac{1}{2}\Im(z\cdot \bar{w})).$
More often we write $ (x,u,t) $ in place of $ (z,t) $ and the group law 
takes the form
$$ (x,u,t)(x',u',t') = (x+x',u+u',t+t'+\frac{1}{2}(u \cdot x'-x \cdot u')) $$
where $ x,u,x',u' \in \R^n.$ For each non-zero
$\lambda \in \R $ the Schrodinger representation $ \pi_\lambda $ of 
$\H^n $ is defined by
$$
 \pi_\lambda(x,u,t)\varphi(\xi) = e^{i\lambda t}e^{i\lambda
(x \cdot \xi+\frac{1}{2}x \cdot u)}\varphi(\xi+u).
$$
The group Fourier transform of $ f \in L^1(\H^n) $ is defined by
$$ \hat{f}(\lambda) = \int_{\H^n} f(z,t)\pi_\lambda(z,t) dz dt.$$
For inversion and palncherel theorems see [13].

As mentioned in the introduction we would like to consider $ \H^n $ as 
the homogeneous space $ G_n/U(n)$  and rewrite the  inversion formula 
in terms of certain representations of $ G_n$. First let us recall some
definitions. The unitary group $ U(n) $ acts on the Heisenberg group as
automorphisms, the action being defined by $ \sigma(z,t) = (\sigma.z,t)$
where $ \sigma \in U(n).$ The Heisenberg motion group $ G_n $ is the
semidirect product of $ U(n ) $ and $ \H^n $ with group law
$$ (\sigma,z,t)(\tau,w,s) = (\sigma \tau,(z,t)(\sigma .w,s)).$$ Functions
on $ \H^n $ can be considered as right $ U(n) $ invariant functions on $ G_n.$
As such the inversion formula for such functions on $ G_n $ will involve only
certain class-one representations of $ G_n.$ We now proceed to describe the
relevant representations.

Let $ \Phi_\alpha, \alpha \in \N^n $ be the normalised Hermite functions
on $ \R^n.$ Let $ \Phi_\alpha^\lambda(x) = |\lambda|^{\frac{n}{4}}\Phi_\alpha
(|\lambda|^{\frac{1}{2}}x) $ and define $ E_{\alpha,\beta}^\lambda(z,t)
= (\pi_\lambda(z,t)\Phi_\alpha^\lambda,\Phi_\beta^\lambda).$ 
For each $ k \in \N $ and non-zero $ \lambda \in \R $ let $ {H}_k^\lambda
$ be the Hilbert space for which the functions $ E_{\alpha,\beta}^\lambda$
with $ \alpha,\beta \in \N^n, |\alpha| =k $ form an orthonormal basis. The inner
product in $ H_k^\lambda $ is defined by
$$ (F,G) = |\lambda|^n \int_{\C^n} F(z,0)\overline{G(z,0)} dz.$$
On this Hilbert space we define a representation $ \rho_k^\lambda $
of the Heisenberg motion group  by
$$
 \rho_k^\lambda(\sigma,z,t)F(w,s) = F((\sigma,z,t)^{-1}(w,s)).
$$
Then it is known that ( see [15])$ \rho_k^\lambda $ is an irreducible unitary
representation of $ G_n.$ As $ (G_n,U(n))$ is a Gelfand pair $ \rho_k^\lambda
$ has a unique $ U(n) $ fixed vector which is none other than 
the Laguerre function $ e_k^\lambda $ (see below).

Given $ f \in L^1(\H^n) $ we can define its group Fourier transform by
$$ \rho_k^\lambda(f) = \int_{G_n} f(z,t)\rho_k^\lambda(\sigma,z,t)
~~d\sigma dz dt $$ which is a bounded operator acting on $ H_k^\lambda \
$.
As shown in [15] we have
$$
 tr(\rho_k^\lambda(\sigma,z,t)^*\rho_k^\lambda(f))
= \frac{k!(n-1)!}{(k+n-1)!} f*e_k^\lambda(z,t)
$$
where $ e_k^\lambda(z,t) = e^{i\lambda t}\varphi_k^\lambda(z)$. Here
$$ \varphi_k^\lambda(z) = L_k^{n-1}(\frac{1}{2}|\lambda||z|^2)e^{-\frac{1}{4}
|z|^2} $$ and $ L_k^{n-1} $ are Laguerre polynomials of type $ (n-1).$ 
The inversion formula for a right $U(n)-$invariant function on $ G_n$
takes the form
$$ f(z,t) = (2\pi)^{-n-1}\int_{-\infty}^\infty \left( \sum_{k=0}^\infty
tr(\rho_k^\lambda(\sigma,z,t)^*\rho_k^\lambda(f))\frac{(k+n-1)!}{k!(n-1)!}
\right) |\lambda|^n d\lambda.$$
Also the Plancherel theorem can be written as
$$
 \int_{\H^n} |f(z,t)|^2 dz dt = \int_{-\infty}^\infty
\left( \sum_{k=0}^\infty
\|\rho_k^\lambda(f)\|_{HS}^2 \frac{(k+n-1)!}{k!(n-1)!}
\right) |\lambda|^n d\mu(\lambda)
$$
where $ d\mu(\lambda) = (2\pi)^{-n-1}|\lambda|^n d\lambda.$

\begin{thm}
For every Schwartz class function $ f $ on $ \H^n $ the following inversion
formula holds:
$$ f(z,t) =\int_{-\infty}^\infty \left( \sum_{k=0}^\infty
\left(\rho_k^\lambda(f)e_k^\lambda,  \rho_k^\lambda(1,z,t)
e_k^\lambda \right) \right) d\mu(\lambda) $$
where $1 $ stands for the identity matrix in $ U(n).$
\end{thm}

From now on let us identify $ \H^n $ with $ \R^n\times\R^n \times \R $
 and use the notation $ (x,u,t) $ rather than $ (x+iu,t)$ to denote
elements of $ \H^n.$ The action of $ U(n) $ on $ \H^n $ then takes the
form $ \sigma.(x,u,t) = (a.x-b.u, b.x+a.u,t) $ where $ a $ and $ b $ are the
real and imaginary parts of $ \sigma.$ This action has a natural extension
to $ \C^n \times \C^n \times \C $ given by 
$\sigma.(z,w,\zeta) = (a.z-b.w, b.z+a.w,\zeta) .$ With this definition we 
can extend the action of
$ G_n $ on $ \H^n $ to $ \C^n \times \C^n \times \C $:
$$ (a+ib,x',u',t')(z,w,\zeta) = (x',u',t')(a.z-b.w, b.z+a.w,\zeta).$$ This action is
then extended to functions defined on  $ \C^n \times \C^n \times \C $:
$$ \rho(g)f(z,w,\zeta) = f(g^{-1}.(z,w,\zeta)),~~~ g \in G_n.$$

We are now ready to prove Gutzmer's formula for the Heisenberg group. Suppose
$ f $ is a Schwartz class function on $ \H^n $ such that $ f^\lambda = 0 $
for all $ |\lambda| > A $ and $ \rho_k^\lambda(f) = 0 $ for all $ \lambda, k $
such that $ (2k+n)|\lambda| > B.$ We say that the Fourier transform of $ f $
is compactly supported if this condition is satisfied for some $ A $ and $ B.$
Now the inversion formula
$$ f(g.(x,u,\xi)) = \int_{-A}^A \sum_{(2k+n)|\lambda|\leq B} (\rho_k^\lambda(f)
e_k^\lambda, \rho_k^\lambda(g)\rho_k^\lambda(1,x,u,\xi)e_k^\lambda) d\mu(\lambda)$$
is valid for any $ g \in G_n.$ Moreover, as each of 
$\rho_k^\lambda(1,x,u,\xi)e_k^\lambda $ extends to $ \C^{2n+1} $ as an entire 
function the same is true of $f(g.(x,u,\xi)) $ and we have
$$f(g.(z,w,\zeta)) = \int_{-A}^A e^{\lambda \eta}\sum_{(2k+n)|\lambda|\leq B} 
(\rho_k^\lambda(f)
e_k^\lambda, \rho_k^\lambda(g)\rho_k^\lambda(1,x,u,\xi)e_k^\lambda) d\mu(\lambda)$$ 
where $ \zeta = \xi+i\eta.$ We then have the following Gutzmer's formula
for the action of Heisenberg motion group on $ \C^{2n+1} $ which is the
complexification of $ \H^n.$

\begin{thm}
Let $ f $ be Schwartz function whose Fourier transform is compactly supported
in the above sense. Then $ f $ extends to $\C^{2n+1} $ as an entire function
and we have the following identity:
$$\int_{G_n} |f(g.(z,w,\zeta))|^2 dg $$ 
$$  =  \int_{-\infty}^\infty 
e^{2\lambda \eta} e^{-\lambda(u\cdot y-v\cdot x)} \left(
\sum_{k=0}^\infty  \|f^\lambda*_\lambda\varphi_k^\lambda\|_2^2
\frac{k!(n-1)!}{(k+n-1)!}
\varphi_k^\lambda(2iy,2iv) \right)d\mu(\lambda) $$
where $\|f^\lambda*_\lambda\varphi_k^\lambda\|_2 $ is the 
$ L^2(\C^n) $ norm of $ f^\lambda*_\lambda\varphi_k^\lambda.$
\end{thm}

For a proof of this theorem we refer to [15] where the formula was proved under
a slightly different condition. In fact, the above formula
holds good as long as the right hand side expression is finite. This can be
proved by means of a density argument.

We now consider the heat kernel transform on the Heisenberg group. Let
$ \CL $ be the sublaplacian on the Heisenberg group and let $ \Delta =
\CL -\partial_t^2 $ be the full Laplacian. Let $ q_t(x,u,\xi) $ be the
heat kernel associated to $ \Delta.$ Then its Fourier transform ( in the
central variable) $ q_t^\lambda(x,u) $ is explicitly known, see [13]. From
the expression it follows that $ q_t(x,u,\xi) $ can be extended to $ \C^{2n+1}
$ as an entire function. The same is true of $ f*q_t(x,u,\xi) $ for any
$ f \in L^2(\H^n).$ In [8] the authors studied the problem of characterising
the image of this heat kernel transform as a space of entire functions on the
complexification $ \H_{\C}^n $ which is just $ \C^{2n+1} .$ They showed
that the image is not a weighted Bergman space but it can be written as
a direct integral of twisted Bergman spaces. They also showed that it is the
direct sum of two Bergman spaces defined in terms if signed weight functions.
Here using Gutzmer's formula we prove another characterisation similar to the
one obtained on Riemannian symmetric spaces. 

Given functions $ m(k,\lambda) $ defined on $ \N \times \R $ we consider
functions of the form
$$ h(iy,iv,i\eta) = \int_{-\infty}^\infty e^{\lambda \eta}\sum_{k=0}^\infty 
m(k,\lambda) \frac{k!(n-1)!}{(k+n-1)!}\varphi_k^\lambda(2iy,2iv) d\mu(\lambda)
.$$ When $ m(k,\lambda) $ is compacktly supported in the sense that it is
supported in $ |\lambda| \leq A, (2k+n)|\lambda| \leq B $ the above
function is well defined and extend to $ \C^{2n+1} $ as an entire function.
Let $ j_{n-1}(s) = s^{-n+1}J_{n-1}(s) $ and define an operator $ \CD $ by
$$ \CD h(iy,iv,i\eta) =
\int_{-\infty}^\infty  e^{\lambda \eta}\sum_{k=0}^\infty    
m(k,\lambda) j_{n-1}(i\sqrt{(2k+n)|\lambda|}(|y|^2+|v|^2)^{\frac{1}{2}})
 d\mu(\lambda)$$
whenever $ h $ is given as above. Note that $ \CD h $ is also an entire 
function on $ \C^{2n+1}.$ Let $ p_t(y,v,\xi) $ be the Euclidean heat
kernel on $ \R^{2n+1}.$ Let us write
$$ O_{|f|^2}(iy,iv,i\eta) = \int_{G_n} |f(g.(iy,iv,i\eta))|^2 dg $$
and call it the orbital integral of $ |f|^2.$ We now state and prove the 
following theorem on the image of the heat kernel transform.

\begin{thm}
An entire function $ F $ on $ \C^{2n+1} $ belongs to the image of the heat
kernel transform on $ L^2(\H^n) $ if and only if
$$ \int_{\R^{2n+1}}\CD O_{|F|^2}(iy,iv,i\eta)p_{t/2}(y,v,\eta)dydvd\eta < \infty.$$ The above is in fact a constant multiple of the $ L^2(\H^n) $ norm of
$ F(x,u,\xi).$
\end{thm} 

{\bf Proof:} Suppose $ F = f*q_t $ for some $ f \in L^2(\H^n).$ Then $ F $
extends to $ \C^{2n+1} $ as an entire function. If $ m(k,\lambda) = \rho_k
^\lambda(F) = e^{-t\lambda^2}e^{-(2k+n)|\lambda|t}\rho_k^\lambda(f) $ then
the function
$$ \int_{-\infty}^\infty e^{2\lambda \eta}\sum_{k=0}^\infty |m(k,\lambda)|^2
 \frac{k!(n-1)!}{(k+n-1)!} \varphi_k^\lambda(2iy,2iv)d\mu(\lambda) $$ 
and hence by Gutzmer's formula $\CD O_{|F|^2}(iy,iv,i\eta) $ is given by
$$ \int_{-\infty}^\infty e^{2\lambda \eta}\sum_{k=0}^\infty  
e^{-2t\lambda^2}e^{-2(2k+n)|\lambda|t} \|f^\lambda *_\lambda \varphi_k^\lambda
\|_2^2 $$
$$ \times
j_{n-1}(2i\sqrt{(2k+n)|\lambda|}(|y|^2+|v|^2)^{\frac{1}{2}}) d\mu(\lambda).
$$
Integrating this against $ p_{t/2}(y,v,\eta) $ and noting that
$$ \int_{\R^{2n+1}} p_{t/2}(y,v,\eta)
j_{n-1}(i\sqrt{(2k+n)|\lambda|}(|y|^2+|v|^2)^{\frac{1}{2}})dy~dv~d\eta$$ 
$$ = e^{2t\lambda^2}
e^{2(2k+n)|\lambda|t} $$ we obtain
$$  \int_{\R^{2n+1}} \CD O_{|F|^2}(iy,iv,i\eta)p_{t/2}(y,v,\eta) dydvd\eta
= c_n \int_{\H^n} |f(x,u,\xi)|^2 dx du d\xi.$$
This proves one half of the theorem. The other half is proved by noting that
all the steps are reversible.

We now state and prove a Paley-Wiener theorem for the inverse Fourier transform
on the Heisenberg group.

\begin{thm}
Let $ f \in L^2(\H^n).$ The Fourier transform $ \rho_k^\lambda(f) $ of $ f $
is compactly supported in $ |\lambda| \leq A, (2k+n)|\lambda| \leq B $ if
and only if $ f $ has an entire extension $ F $ to $ \C^{2n+1} $ which
satisfies the estimate
$$ \CD O_{|F|^2}(iy,iv,i\eta) \leq C e^{2A|\eta|}e^{2\sqrt{B}(|y|^2+|v|^2)^
{\frac{1}{2}}} $$
for all $ (y,v,\eta) \in \R^{2n+1}.$
\end{thm}

{\bf Proof:} First assume that $ \rho_k^\lambda(f) $ 
is compactly supported in $ |\lambda| \leq A, (2k+n)|\lambda| \leq B.$ As we
have seen in the proof of Gutzmer's formula $ f $ extends to an entire function
$ F $  and $ \CD O_{|F|^2}(iy,iv,i\eta) $ is given by
$$\int_{-A}^A e^{2\lambda \eta}\sum_{(2k+n)|\lambda|\leq B}
\|f^\lambda *_\lambda \varphi_k^\lambda \|_2^2
j_{n-1}(2i\sqrt{(2k+n)|\lambda|}(|y|^2+|v|^2)^{\frac{1}{2}}) d\mu(\lambda). $$
As $ j_{n-1}(is)\leq Ce^s $ the above gives the estimate
$$ \CD O_{|F|^2}(iy,iv,i\eta) \leq C e^{2A|\eta|}e^{2\sqrt{B}(|y|^2+|v|^2)^
{\frac{1}{2}}} \|f\|_2^2.$$ This proves the sufficiency part of the theorem.

To prove the necessity, assume that $ F $ satisfies the hypothesis of the
theorem. First of all the Euclidean Paley-Wiener theorem for the central
variable shows that $ f^\lambda = 0 $ for all $ |\lambda| > A.$  The
hypothesis then implies that $ F $ belongs to the image of the heat kernel
transform $ f \rightarrow f*q_t $ for any $ t > 0 $ and also
$$ \int_{\R^{2n+1}}\CD O_{|F|^2}(iy,iv,i\eta)p_{t/2}(y,v,0) p_{s/2}(0,0,\eta)
dy dv d\eta \leq C e^{2sA^2} e^{2tB}.$$ By Gutzmer's formula this means that
$$ \int_{-\infty}^\infty e^{2s\lambda^2} \sum_{k=0}^\infty \|f^\lambda *_
\lambda \varphi_k^\lambda\|_2^2 e^{2(2k+n)|\lambda|t} d\mu(\lambda)
\leq C e^{2sA^2} e^{2tB}.$$ As this is true for every $ t $ proceeding as in
the case of symmetric spaces we can show that $\rho_k^\lambda(f) $ is
supported in $ (2k+n)|\lambda| \leq B.$ This completes the proof of the
theorem.

The following  remarks on the operator $ \CD $ are in order. It has the 
effect of changing $ \varphi_k^\lambda(iy,iv) $ into 
$ j_{n-1}(i\sqrt{(2k+n)|\lambda|}(|y|^2+|v|^2)^{\frac{1}{2}})$ Notice that
$\varphi_k^\lambda $ are the spherical functions associated to the $k$-the ray
of the Heisenberg fan whereas $ j_{n-1} $ is the spherical function
associated to the limiting ray. Moreover, from the asymptotic formula
of Hilb's type for Laguerre functions ( see Theorem 8.22.4 in Szego [11] )
we see that $\varphi_k^\lambda(iy,iv) $ is approximated by
$ j_{n-1}(i\sqrt{(2k+n)|\lambda|}(|y|^2+|v|^2)^{\frac{1}{2}})$. Further study
of the operator $ \CD $ is worth considering. 

In view of the above remarks the above theorem is not completely satisfactory.
We can prove another version of Paley-Wiener theorem if we make use of the
characterisation of the image of $ L^2(\H^n) $ under the heat kernel transform
obtained in [8]. There the authors have shown that the image is the direct sum
of two Bergman spaces each of which is defined interms of certain weight 
function which takes both positive and negative values. To be more precise
these weight functions denoted by $ W_t^+ $ and $ W_t^- $ are defined by the
equations
$$ \int_{-\infty}^\infty e^{2\lambda \eta} W_t^+(iy,iv,\eta) d\eta =
e^{2t\lambda^2}p_{2t}^\lambda(2y,2v) $$
for $ \lambda >0,$ and a similar equation for $  W_t^- $ valid for 
$ \lambda < 0.$ In the above $ p_t^\lambda $ is the heat kernel associated to
to the special Hermite operator and given explicitly by
$$ p_t^\lambda(z,w) = c_n \left( \frac{\lambda}{\sinh (\lambda t)}\right)^n
e^{-\frac{\lambda}{4}\coth (\lambda t)(z^2+w^2)}.$$ The existence of such
weight functions have been proved in [8].

We now have the following theorem. We consider only functions $ f $ for which $
f^\lambda $ is supported in $ \lambda > 0.$  A similar result is true for functions $ f $ for which $ f^\lambda $ is supported in $ \lambda < 0.$

\begin{thm} Let $ f \in L^2(\H^n) $ be such that $ f^\lambda $ is supported in
 $ \lambda > 0.$ If $ \rho_k^\lambda(f) $ is compactly supported then $ f $
extends to $ \C^{2n+1} $ as an entire function $ F(z,w,\zeta) $ which is of
exponential type in the last variable and satisfies the following condition:
for some constants $ B, C >0 $ we have
$$ |\int_{-\infty}^\infty O_{|F|^2}(iy,iv,i\eta)W_t^+(2iy,2iv,i\eta)d\eta|
\leq C t^{-2n}e^{2tB} $$
for all $ t> 0.$ Conversely, if $ F(z,w,\zeta) $ is entire, of exponential type
in $ \zeta $ and satisfies the slightly stronger estimate
$$ |\int_{-\infty}^\infty O_{|F|^2}(iy,iv,i\eta)W_t^+(iy,iv,i\eta)d\eta|
\leq C_N t^{-2n}(1+|y|^2+|v|^2)^{-N}e^{2tB} $$
for some $ N > n $ and for all $ t > 0 $ then $ \rho_k^\lambda(f) $ is
compactly supported.
\end{thm}

{\bf Proof:} If $ \rho_k^\lambda(f) $ is supported in $ 0 < \lambda \leq \alpha
$ and $ (2k+n)|\lambda| \leq \beta $ then by Gutzmers' formula and the defining
relation for $ W_t^+ $ we have
$$\int_{-\infty}^\infty O_{|F|^2}(iy,iv,i\eta)W_t^+(2iy,2iv,i\eta)d\eta $$
$$ = \int_{0}^\alpha e^{2t\lambda^2} \sum_{(2k+n)|\lambda| \leq \beta} \|f^\lambda
*_\lambda \varphi_k^\lambda\|_2^2 \frac{k!(n-1)!}{(k+n-1)!}
\varphi_k^\lambda(2iy,2iv)p_{2t}^\lambda(
4y,4v) d\mu(\lambda).$$  
Now the function $\varphi_k^\lambda(z,w) $ belongs to the twisted Bergman
spaces $ \CB_t^\lambda $ studied in [8] for any $ t > 0.$ The reproducing
kernel $ K_t^\lambda $ for these spaces are given in terms of $ p_{2t}^\lambda $ (see [8]): more precisely
$$ K_t^\lambda((z,w),(a,b)) = p_{2t}^\lambda(z-\bar{a},w-\bar{b})e^{-\frac{i}{2}\lambda (w \cdot \bar{a} - z \cdot \bar{b})}.$$
As evaluations are continuous it follows that
$$ |\varphi_k^\lambda(2iy,2iv)| \leq K_t^\lambda((2iy,2iv),(2iy,2iv))
\|\varphi_k^\lambda \| .$$ Since $ \|\varphi_k^\lambda\| = C\frac{(k+n-1)!}
{k!(n-1)!}
 e^{2(2k+n)|\lambda|t}$ it follows that
$$ \frac{k!(n-1)!}{(k+n-1)!} \varphi_k^\lambda(2iy,2iv)p_{2t}^\lambda(4y,4v) 
 \leq C \left(\frac{\lambda}{\sinh(2t\lambda)}\right)^{2n}
e^{2(2k+n)|\lambda|t}.$$
Using this estimate in the above  and noting that $ \frac{\lambda}
{\sinh \lambda} $ is a decreasing function we obtain
$$ \int_{-\infty}^\infty O_{|F|^2}(iy,iv,i\eta)W_t^+(2iy,2iv,i\eta)d\eta $$
$$\leq C t^{-2n} e^{t \alpha^2} e^{2t\beta} \|f\|_2^2.$$ This proves one half of the theorem with $ B = \alpha^2 + \beta.$

Now for the converse. The hypothesis on $ F $ means that
$$ \int_{0}^\infty  e^{2t\lambda^2} \sum_{k=0}^\infty 
 \|f^\lambda *_\lambda \varphi_k^\lambda\|_2^2 \frac{k!(n-1)!}{(k+n-1)!}
\varphi_k^\lambda(2iy,2iv)p_{2t}^\lambda(
2y,2v) d\mu(\lambda)$$
$$ \leq  C_N t^{-2n}(1+|y|^2+|v|^2)^{-N}e^{tB}.$$ Integrating with respect to
$ dy dv $ and noting that
$$ \int_{\R^{2n}} \varphi_k^\lambda(iy,iv)p_t^\lambda(y,v) dy dv =
c_n \frac{(k+n-1)!}{k!(n-1)!} e^{(2k+n)|\lambda|t} $$ 
(see Lemma 6.3 in [15]) we have
$$ \int_{0}^\infty  e^{2t\lambda^2} \sum_{k=0}^\infty
 \|f^\lambda *_\lambda \varphi_k^\lambda\|_2^2 e^{2(2k+n)|\lambda|t} 
d\mu(\lambda) 
 \leq  C_N t^{-2n}e^{2tB}.$$ This gives for ant $ C > 0 $ 
$$ e^{2tC}\int_{0}^\infty  \sum_{(2k+n)|\lambda|> C}
 \|f^\lambda *_\lambda \varphi_k^\lambda\|_2^2 
d\mu(\lambda)  \leq C_N e^{2tB} $$
for all $ t \geq 1.$ Taking $ C > B $ and letting $ t $ tend to infinity we
conclude that $ \rho_k^\lambda(f) $ is supported in $ (2k+n)|\lambda| \leq B.
$ This completes the proof.

In the above theorem the necessary and sufficient conditions are different. For
the necessary condition involves pointwise estimate on the spherical functions
$ \varphi_k^\lambda $ whereas for the sufficiency we have used an integral
condition on the same functions. We therefore, may not hope to get the same
condition as both necessary and sufficient.

\end{document}